\documentclass{osa-article}

\journal{oe}

\usepackage[symbol]{footmisc}

\newtheorem{theorem}{Theorem}

\newtheorem{proposition}[theorem]{Proposition}

\articletype{Research Article}

\begin{document}

\title{Hyperspectral holography and spectroscopy: computational features of inverse
discrete cosine transform }

\author{Vladimir Katkovnik, Igor Shevkunov and Karen Egiazarian}
\address{Faculty of Information Technology and Communication, Tampere University, Tampere, Finland}

\email{vladimir.katkovnik@tuni.fi} 



\begin{abstract}

Broadband hyperspectral digital holography and Fourier transform spectroscopy
are important instruments in various science and application fields. In the
digital hyperspectral holography and spectroscopy the variable of interest are obtained as
inverse discrete cosine transforms of observed diffractive intensity patterns.
In these notes, we provide a variety of algorithms for the inverse cosine
transform with the proofs of perfect spectrum reconstruction, as well as
we discuss and illustrate some nontrivial features of these algorithms.
\end{abstract}
\section{Introduction}

In the hyperspectral complex domain (phase/amplitude) imaging based on
phase-shifting holography there are two testing light beams, object $A$ and
reference $R$, with the object placed in the object beam. Hologram
(interferogram) patterns with varying distance (phase) shifts between the
object and reference beams arriving to the sensor are registered. These
observation are used for the object reconstruction provided that the reference
beam is known (\cite{Juptner-2005}, \cite{Review-2010}). For the broadband
object and reference beams (white light) the registered holograms have the
following form (\cite{Kalenkov-2012}, \cite{Kalenkov-2017}):
\begin{align}
&  J(t)=\int_{0}^{\infty}|A(\omega)+R\exp(j2\pi\omega t)|^{2}d\omega
=\label{Integral-Interferogram}\\
&  \int_{0}^{\infty}(|A(\omega)|^{2}+R^{2}+2|A(\omega)|\times R\cos
(\varphi_{A}(\omega)-2\pi\omega t))d\omega,\nonumber
\end{align}

where $A(\omega)$ is a complex-valued transfer function of the object (for
transparent object) to be analyzed with the amplitude $|A(\omega)|$ and the
phase $\varphi_{A}(\omega)$, $A(\omega)=|A(\omega)|\exp(j\varphi_{A}(\omega))$.

The similar model is used for a reflective object where $A(\omega)$ is a
complex-valued spectrum of the reflective light beam. $R$ is the amplitude of
the reference beam invariant with respect to $\omega$ and $\varphi
(\omega)=2\pi\omega t$ is varying phase-shift of the reference beam with
respect to the object beam.

This phase-shift is controlled by the parameter $t$. The $\omega=1/\lambda$
stays for the frequency corresponding to the wavelength $\lambda$. Integration
over $\omega$ is produced over the spectrum of the broadband object beam
defined by the spectrum of the input light beam impinging on the object.

The integration on the non-negative $\omega$ in
Eq.(\ref{Integral-Interferogram}) takes arguments in physics while in
mathematics the symmetric representation is more usual with integration from
all frequencies from minus to plus infinity.

Here we use the non-symmetric version of this representation corresponding to
(\ref{Integral-Interferogram}).

The optical setup conventional for spectroscopy without a reference beam
assumes that both the basic and shifted broadband beams go through the object.
Thus, two identical but mutually shifted broadband copies of the object are
superposed on the sensor plane ( \cite{Bell-1972}, \cite{2001-Davis-Book}). In
\cite{Falldorf}, this lateral shift is referred to as the shear. Citing
\cite{Falldorf}, the \textquotedblright s\textit{hear interferometry has a
great advantage over the conventional one based on a reference wave, because
the interfering wave fields travel almost the same path and as a result the
shear interferometry is insensitive to environmental disturbances and has very
low demands with respect to the coherence of the light}.\textquotedblright

The observed intensity patterns usually are presented in the form
\cite{2001-Davis-Book}:
\begin{align}
&  J(t)=\int_{0}^{\infty}|A(\omega)+A(\omega)\exp(j2\pi\omega t)|^{2}%
d\omega=\label{Spectrometry-Integral}\\
&  2\int_{0}^{\infty}|A(\omega)|^{2}(1+\cos(2\pi\omega t))d\omega,\nonumber
\end{align}

where again $\varphi(\omega)=2\pi\omega t$ is the varying phase-shift
controlled by the parameter $t$.

In the both equations (\ref{Integral-Interferogram}) and
(\ref{Spectrometry-Integral}) the observations $J(t)$ are the cosine
transforms of the variables of interest $A(\omega)$ for
(\ref{Integral-Interferogram}) and $|A(\omega)|^{2}$ for
(\ref{Spectrometry-Integral}). The cosine transform calculated for $J(t)$ over
$\omega$ with the shift variable $t$, for instance as $\int J(t)\cos
(2\pi\omega t)dt$, is the standard instrument in order to extract spectral
information on the object $A(\omega) $ from the observed $J(t)$.

Mainly, the corresponding algorithms are presented and discussed for the
integral forms of the observations shown above (e.g. \cite{Kalenkov-2017}),
while the practical calculations are produced using the Fast Fourier Transform
(FFT) algorithm, which is one of the versions of the discrete Fourier transforms (DFTs) which are quite different
from the integral cosine and Fourier transforms.

In these notes, we are focused on specific features of FFT applied to the
discrete versions of (\ref{Integral-Interferogram}) and
(\ref{Spectrometry-Integral}) appeared due to the periodical nature of FFT and
observations $J(t).$

Our goal is to derive the precise estimates of spectra from given discrete observations.
Thus, as the first step we replace (\ref{Integral-Interferogram}%
)-(\ref{Spectrometry-Integral}) by the corresponding discrete models
corresponding to the conventions of FFT and derive the algorithms for precise
reconstruction of the corresponding discretized spectra.

\section{Spectroscopy: reference less setup}

Let us consider instead of spectrum reconstruction for
(\ref{Spectrometry-Integral}) a more general problem: inverse of the cosine
transform formalized as
\begin{equation}
J(t)=2\int_{0}^{\infty}X(\omega)(1+\cos(2\pi\omega t))d\omega,
\label{Spectrometry-Integral1}%
\end{equation}

where the unknown $X(\omega)$ is real-valued and may take both positive and
negative values.

Keeping in mind that the estimates will be applied for $X(\omega
)=|A(\omega)|^{2}$, we continue to use for $X(\omega)$ the term
\textit{spectrum}. The inverse cosine transform problem is to find $X(\omega)$
from observations Eq.(\ref{Spectrometry-Integral1}).

\subsection{Non-symmetric sampling\label{Nonsymmetric sampling}}

Following the convention of FFT, assume that the
sampling intervals $\Delta_{t}$ for the shift  and $\Delta_{\omega}$ for the frequency
 are such that
\begin{equation}
N=\frac{1}{\Delta_{t}\Delta_{\omega}}, \label{Nauquist}%
\end{equation}
with a number of samples $N$ for both $J(t)$ and $X(\omega)$.

The discretization of (\ref{Spectrometry-Integral1}) is produced according to
the relations $t=\Delta_{t}\tau$ and $\omega=\Delta_{\omega}u$, where $\tau$
and $u$ are integers. Integration for $J(t)$ in (\ref{Spectrometry-Integral1})
is done over the finite interval $[0,\Delta_{\omega}(N-1)]$ with the sampled
values of $J(\tau)$ over the interval of integers $\tau\in\lbrack1,2,...,N].$

We call this sampling as non-symmetric contrary to the symmetric one with \\
$\tau\in\lbrack-N/2,...,N/2-1]$ considered in the next subsection.

The discrete version of the formula (\ref{Spectrometry-Integral}) takes the
form%
\begin{equation}
J(\tau)=2\sum_{u=0}^{N/2-1}X(u)(1+\cos(\frac{2\pi}{N}\tau u)),\text{ }%
\tau=1,..N, \label{SpectrallyResolvedObs}%
\end{equation}
where the integer discrete frequency $u$ covers the low frequency interval
$\{0,1,...,N/2-1\}$.

The restriction of the frequency interval for the estimated spectrum to the
length $N/2$ is conventional to work with FFT and follows from the periodicity
of FFT, see the proof of Proposition 1 in Appendix.

Thus, the signal to be spectrally analyzed from observations $J(\tau)$ should
be \textit{'low band'} of the length $N/2$ with $A(u)=0$ for $N/2\leq u\leq
N-1$ provided that the sampling length is $N$.

We say that the observations $J(\tau)$ are spectrally resolved if $X(u)$,
$0\leq u\leq N/2-1$, satisfying to Eq.(\ref{SpectrallyResolvedObs}) are found.

\begin{proposition}
Let the observations $J(\tau)$, $\tau=1,...,N$, be defined as in
(\ref{SpectrallyResolvedObs}), then $X(u)$ can be calculated as follows:%
\begin{equation}
X(u)=\left\{
\begin{array}
[c]{c}%
\dfrac{1}{N}\sum_{\tau=1}^{N}J(\tau)\cos(\frac{2\pi}{N}\tau u)\text{, if
}u=1,..,N/2-1,\\
\\
\dfrac{1}{4N}\sum_{\tau=1}^{N}J(\tau)-\frac{1}{2}\sum_{u_{1}=1}^{N/2-1}%
X(u_{1}),\text{ if }u=0.
\end{array}
\right.  \label{A_tay_1}%
\end{equation}

Let FFT be used for calculations and
\begin{align}
x(u)  &  =\operatorname{real}(IFFT(J(\tau))\times\exp(j\frac{2\pi}%
{N}(u-1))),\label{a1}\\
u  &  =1,..,N,\nonumber
\end{align}
where IFFT stands for the inverse FFT and the variable $u$ is the  argument of $IFFT(J(\tau))$.

Then, the formula (\ref{A_tay_1}) can be written as
\begin{equation}
X(u)=\left\{
\begin{array}
[c]{c}%
x(u+1)\text{, if }u=1,...,N/2-1,\\
\\
\dfrac{1}{4N}\sum_{\tau=1}^{N}J(\tau)-\frac{1}{2}\sum_{u_{1}=1}^{N/2-1}%
X(u_{1}),\text{ if }u=0.
\end{array}
\right.  \label{A_tay_1_IFFT}%
\end{equation}

\end{proposition}

The proof of the formulas (\ref{A_tay_1})-(\ref{A_tay_1_IFFT}) is given in Appendix.

The MATLAB convention for FFT and
IFFT is assumed here and in what follows:
\begin{align}
FFT(y_{\tau})  &  :Y_{u}=\sum_{\tau=1}^{N}y_{\tau}\exp(-j\frac{2\pi}{N}%
(\tau-1)(u-1)),\label{FFT_IFFT}\\
u  &  =1,..,N;\nonumber\\
IFFT(Y_{u})  &  :y_{\tau}=\frac{1}{N}\sum_{\tau=1}^{N}Y_{u}\exp(j\frac{2\pi
}{N}(\tau-1)(u-1)),\nonumber\\
\tau &  =1,..,N.\nonumber
\end{align}

The exponent $\exp(j\frac{2\pi}{N}(u-1))$ in (\ref{a1}) compensates the
difference between this definition of IFFT and the arguments we need for (\ref{A_tay_1}).

The formula (\ref{A_tay_1_IFFT}) looks unusually being different by the
expression for $X(0)$ and by the exponential correction for IFFT. Computational tests demonstrate that if we omit the exponential
correction the results degrade and do not correspond to $J(\tau)$ defined by
(\ref{SpectrallyResolvedObs}). Even more, if we assume that $X(u)\neq0$ for
$n>N/2$, the estimates for $X(u)$ degrade quickly for $u\leq N/2-1$ as
soon as $X(u)\neq0$ for larger $u.$

In what follows, we discuss a few simulation tests illustrating these specific
features of the derived estimates of the spectrum. The estimates are
calculated using IFFT as in (\ref{a1})-(\ref{A_tay_1_IFFT}).

\begin{figure}[ht]%
\centering
\includegraphics[
height=2.5763in,
width=4.5071in
]%
{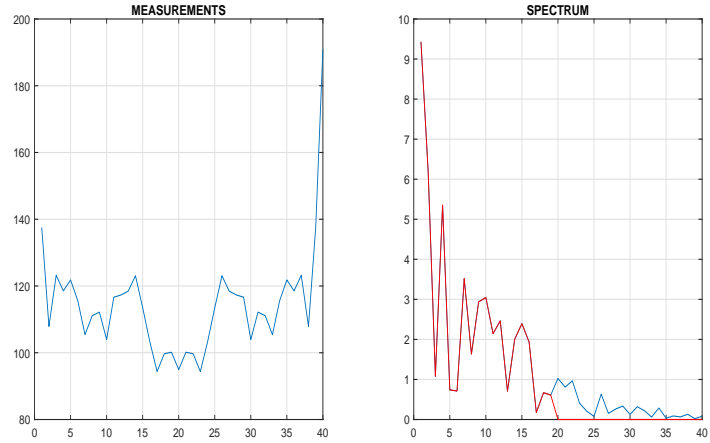}%
\caption{True spectrum $X(u)$ and the corresponding observations $J(\tau)$,
$N=40$.}%
\label{spectr_illustr_11}%
\end{figure}

\begin{figure}[ht]%
\centering
\includegraphics[
height=2.5763in,
width=4.5071in
]%
{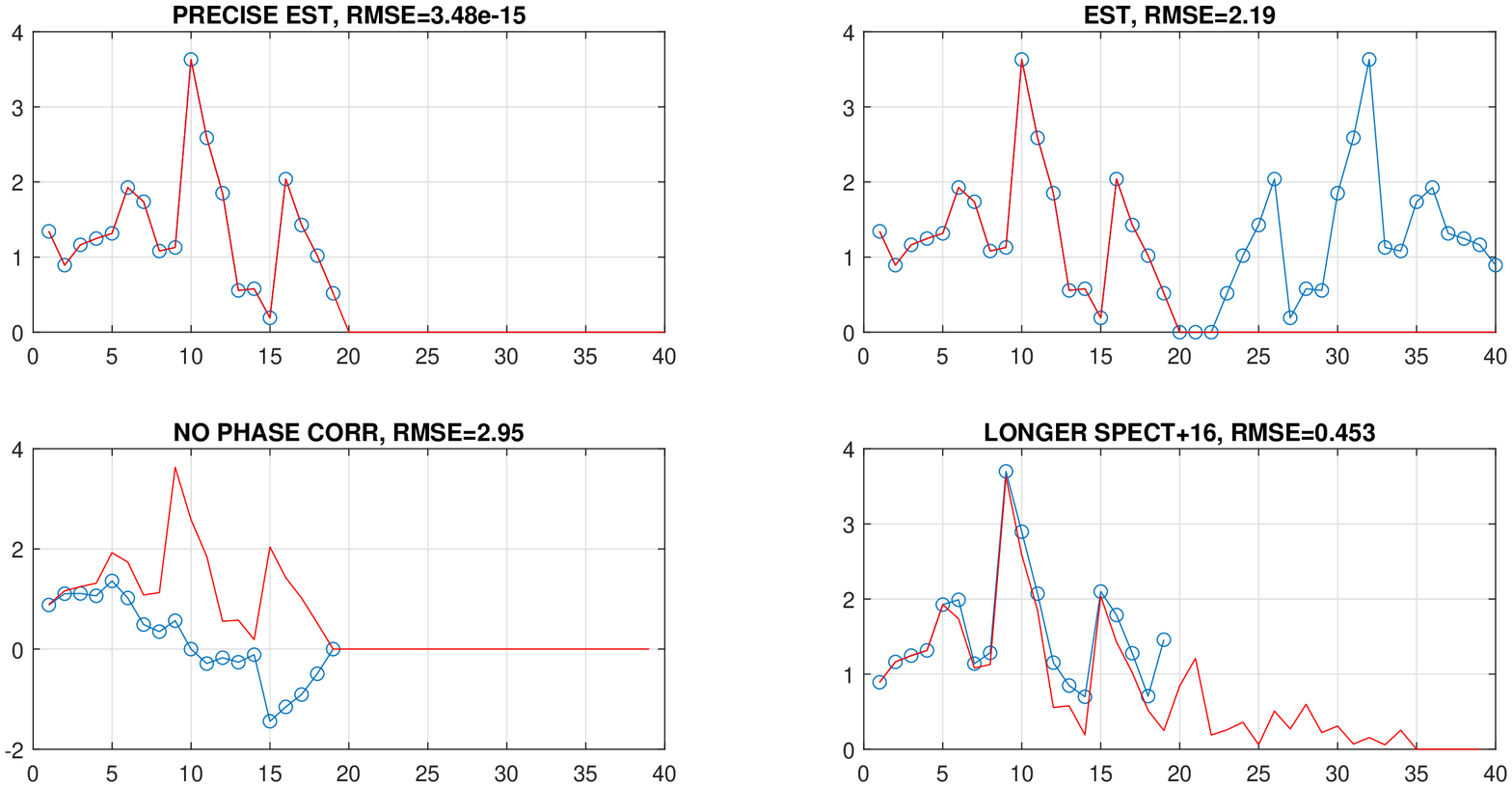}%
\caption{Spectrum reconstructions: the true values are in red and the
estinates are blue with circles, $N=40.$}%
\label{estimate_illustr_1}%
\end{figure}

In Fig.\ref{estimate_illustr_1}, we show the estimates obtained for the data
in Fig.\ref{spectr_illustr_11}, where the measurements are given by the
variable $J(t)$ and the spectrum corresponds to $X(u)$, the number of
observations $N=40$. According to the above results the perfect reconstruction
is achieved for for $0\leq u\leq N/2-1$ provided that $X(u)=0$ for $u\geq
N/2$. The corresponding observations are obtained from $X(u)$ shown in Fig.\ref{spectr_illustr_11} by zeroing values of $X(u)$ for $u\geq N/2$. This
truncated $X(u)$ is marked by red in Fig.\ref{spectr_illustr_11} and the part
of the spectrum $X(u)$ for $u\geq N/2$ is blue.%

The precise estimates shown in Fig. \ref{estimate_illustr_1} are obtained for
this truncated spectrum zeroed for $u\geq N/2$ and shown in the upper row of
the images. In the left image of this row, we show the estimates calculated for
the basic interval $0\leq u\leq N/2-1.$ The perfect reconstruction is
demonstrated in this case with the very small value of RMSE. In images, we show
the true values and the estimates for $u\geq1$ as the values for $u=0$ have
large absolute values while RMSE are calculated for the interval $u=0,..,N/2-1$.

In the right image of the upper row, we show the estimates calculated for all $0\leq u\leq N-1
$. The red color shows the true values of the spectrum, thus we can see that
indeed the estimates calculated for $u\geq N/2$ are completely erroneous
because the true values of the used spectra are equal to zero. There is an
obvious symmetry of the estimates for $u\geq N/2$ with respect to the
estimates for $u<N/2$. The right part of the estimate is obtain by mirror
reflection with respect to the point $N/2+1.$%

The second row in Fig.\ref{estimate_illustr_1} illustrates two different effects. The left estimate is calculated as it is in the first row
but without the phase correction for $IFFT$ in Eq.(\ref{a1}). One may note
that this phase correction is of importance as the estimate is completely
destroyed, what can be noticed from the large values of RMSE as well as
visually comparing the estimate and the true spectral values.

The last image shows the influence of the requirement $X(u)=0$ for $u\geq
N/2.$ Here we assume that $X(u)=0$ for $u\geq N/2+16$. Thus, $16$ non-zero
pixels of this enlarged $X(u)$ influence the $40$ observations of $J$, which are
naturally different from those in Fig.\ref{spectr_illustr_11}. The estimates
and the corresponding true $X(u)$ are seen in Fig. \ref{estimate_illustr_1}. We may
note that the estimates are destroyed with RMSE much larger than those for the
perfect reconstruction demonstrated in the first row of the images.

\bigskip

Let us consider a set of observations $J(\tau)$ in the form%
\begin{equation}
J(\tau)=2\sum_{u=0}^{N/2-1}X(u)(1+\cos(\frac{2\pi}{N}\tau u)),\text{ }%
\tau=0,1,..N-1, \label{SpectrallyResolvedObs2}%
\end{equation}
which is \textit{different\ from} (\ref{SpectrallyResolvedObs}) by the set of
shifted $\tau$ including now $\tau=0$.

This $\tau=0$ results in significant changes in the reconstruction algorithm implemented using FFT.

\begin{proposition}
Let the observations $J(\tau)$, $\tau=0,1,..N-1$, be defined as in
(\ref{SpectrallyResolvedObs2}), then the spectra $X(u)$ can be calculated as
follows:%
\begin{equation}
X(u)=\left\{
\begin{array}
[c]{c}%
\dfrac{1}{N}\sum_{\tau=0}^{N-1}J(\tau)\cos(\frac{2\pi}{N}\tau u)\text{, if
}u=1,..,N/2-1,\\
\\
\dfrac{1}{4N}\sum_{\tau=1}^{N}J(\tau)-\frac{1}{2}\sum_{u_{1}=1}^{N/2-1}%
X(u_{1}),\text{ if }u=0.
\end{array}
\right.  \label{A_tay_12}%
\end{equation}

Let FFT be used for calculations of (\ref{A_tay_12}) and
\begin{equation}
x(u)=\operatorname{real}(IFFT(J(\tau))),\text{ }u\in\lbrack1,...,N],
\label{a12}%
\end{equation}
then
\begin{equation}
X(u)=\left\{
\begin{array}
[c]{c}%
x(u+1)\text{, if }u=1,..,N/2-1,\\
\\
\dfrac{1}{4N}\sum_{\tau=0}^{N-1}J(\tau)-\frac{1}{2}\sum_{u_{1}=1}%
^{N/2-1}X(u_{1}),\text{ if }u=0.
\end{array}
\right.  \label{A_tay_1_IFFT2}%
\end{equation}

\end{proposition}

The proofs of the formulas (\ref{A_tay_12})-(\ref{A_tay_1_IFFT2}) are similar
to given for the proof of Proposition 1 in Appendix.

Thus, summation on $\tau$ in (\ref{A_tay_12}) is produced from $\tau=0$ up to
$\tau=N-1$, what is the only difference  from (\ref{A_tay_1}). The more serious
difference appeared in (\ref{a12}), where the phase correction specific for
(\ref{a1}) is not required.

The solutions given by the formulas (\ref{A_tay_1})-(\ref{A_tay_1_IFFT}) and
(\ref{A_tay_12})-(\ref{A_tay_1_IFFT2}) are precise solutions of the equations
(\ref{SpectrallyResolvedObs}) and (\ref{SpectrallyResolvedObs2}) linear with
respect to $X(u)$.

The properties of the estimates (\ref{A_tay_12})-(\ref{A_tay_1_IFFT2}) are
identical to those discussed above for solutions (\ref{A_tay_1}%
)-(\ref{A_tay_1_IFFT}) excluding the phase correction in (\ref{a1}) which does
not exist in (\ref{a12}).

\begin{proposition}
Let in both (\ref{SpectrallyResolvedObs}) and (\ref{SpectrallyResolvedObs2})
$X(u)$ be real-valued non-negative, i.e. $X(u)=|A(u)|^{2}$, the estimates
(\ref{A_tay_1}) and (\ref{A_tay_12}) continue to be valid in these cases\ and
the estimates (\ref{A_tay_1_IFFT}) and (\ref{A_tay_1_IFFT2}) given by FFT are
simplified to the following form valid for the both estimates:%
\begin{equation}
a(u)=abs(IFFT(J(\tau))),\text{ }u\in\lbrack1,...,N], \label{a12-final}%
\end{equation}%
\begin{equation}
|A(u)|^{2}=\left\{
\begin{array}
[c]{c}%
a(u+1)\text{, if }u=1,..,N/2-1,\\
\\
\dfrac{1}{4N}\sum_{\tau=0}^{N-1}J(\tau)-\frac{1}{2}\sum_{u_{1}=1}%
^{N/2-1}|A(u_{1})|^{2},\text{ if }u=0.
\end{array}
\right.  \label{A_tay_1_IFFT2-final}%
\end{equation}

\end{proposition}

The estimates (\ref{A_tay_1}) and (\ref{A_tay_12}) are valid for
$X(u)=|A(u)|^{2}$ just because it is a particular case for $X(u)$. In
(\ref{a12-final}) the operation $real$ is replaced by $abs$ just because
$|A(u)|^{2}$ are non-negative. Then, the complex exponent in (\ref{a1}) is
killed by $abs$ and the estimates in both Proposition 1 and Proposition 2 take
the identical form presented above in Proposition 3.

Note, that $a(u)$ in (\ref{a12-final}) can be calculated using $FFT$ as 
\begin{equation}
a(u)=abs(FFT(J(\tau)))/N,\text{ }u\in\lbrack1,...,N]. \label{a12-final11}%
\end{equation}%

\subsection{Symmetric sampling\label{Symmetric sampling}}

Now, we consider the symmetric sampling of $J(\tau)$ with $\tau\in
\lbrack-N/2,...,N/2+1]$.

Then the observation model (\ref{SpectrallyResolvedObs}) takes the form%
\begin{equation}
J(\tau)=2\sum_{u=0}^{N/2-1}X(u)(1+\cos(\frac{2\pi}{N}\tau u)),\text{ }%
\tau=-N/2,..N/2-1. \label{sym-sampling}%
\end{equation}

The signal to be spectrally analyzed by observations $J(\tau)$ is of 'low
band' as it is in (\ref{SpectrallyResolvedObs}).

We say that the observations $J(\tau)$ are spectrally resolved if $X(u) $
satisfying to Eq.(\ref{sym-sampling}) are found.

\begin{proposition}
Let the observations $J(\tau)$, $\tau=-N/2,..N/2-1$, be defined according to
(\ref{sym-sampling}), then $X(u)$ can be calculated as follows:%
\begin{equation}
X(u)=\left\{
\begin{array}
[c]{c}%
\frac{1}{N}\sum_{\tau=-N/2}^{N/2-1}J(\tau)\cos(\frac{2\pi}{N}\tau u)\text{, if
}u=1,..,N/2-1,\\
\\
\frac{1}{4N}\sum_{\tau=-N/2}^{N/2-1}J(\tau)-\frac{1}{2}\sum_{u_{1}=1}%
^{N/2-1}X(u_{1}),\text{ if }u=0.
\end{array}
\right.  \label{sym-spectra-est2}%
\end{equation}
If FFT is used for calculations these estimates can be presented in the following two forms:

(1) Let
\begin{equation}
x(u)=\operatorname{real}(IFFT(J(\tau))\times\exp(-j\pi(u-1))), u=1,...,N,
\label{x_u_old}%
\end{equation}
then
\begin{equation}
X(u)=\left\{
\begin{array}
[c]{c}%
x(u+1)\text{, }u=1,2,...,N/2-1\text{,}\\
\\
\frac{1}{4N}\sum_{\tau=-N/2}^{N/2-1}J(\tau)-\frac{1}{2}\sum_{u=1}%
^{N/2-1}|A(u)|^{2},\text{ if }u=0;
\end{array}
\right.  \label{sym-fft-old}%
\end{equation}
\ 

(2) Let
\begin{equation}
x(u)=\operatorname{real}(IFFT(fftshift(J(\tau)))),\text{ }u\in\lbrack
1,...,N], \label{a22}%
\end{equation}
then
\begin{equation}
X(u)=\left\{
\begin{array}
[c]{c}%
x(u+1)\text{, }u=1,2,...,N/2-1\text{,}\\
\\
\frac{1}{4N}\sum_{\tau=-N/2}^{N/2-1}J(\tau)-\frac{1}{2}\sum_{u_{1}=1}%
^{N/2-1}X(u_{1}),\text{ if }u=0.
\end{array}
\right.  \label{sym-fft2}%
\end{equation}

\end{proposition}
The proof of the proposition can be seen in Appendix.

Note, that the phase correction in (\ref{x_u_old}) is essential as it changes
signs of the items of IFFT, $\exp(-j\pi(u-1)))$ periodically takes values $1$
and $-1$, respectively for $u=1,...,N.$

This kind of the phase correction is not required for (\ref{a22}%
)-(\ref{sym-fft2}).

The shifts on $\Delta\tau$ of the argument\ in $J(\tau)$ results in the
complex exponent factor for the corresponding FFT, $FFT(J(\tau+\Delta
\tau))=FFT(J(\tau))\times\exp(-j\frac{2\pi}{N}\Delta\tau u)$ and
$abs(FFT(J(\tau+\Delta\tau)))=abs(FFT(J(\tau)))$.

It explains, that for the problems in Propositions 1,2,3 we have identical
solutions when we talk about the reconstruction of the spectrum $|A(u)|^{2}$
through FFT.

With the same reasons, the estimates (\ref{a22})-(\ref{sym-fft2}) for
$X(u)=|A(u)|^{2}$ can be replaced with
\begin{equation}
a(u)=abs(IFFT(\textit{fftshift}(J(\tau))))),\text{ }u=1,...,N. \label{a333}%
\end{equation}
and
\begin{equation}
|A(u)|^{2}=\left\{
\begin{array}
[c]{c}%
a(u+1)\text{, }u=1,2,...,N/2-1\text{,}\\
\\
\frac{1}{4N}\sum_{\tau=-N/2}^{N/2-1}J(\tau)-\frac{1}{2}\sum_{u_{1}=1}%
^{N/2-1}|A(u_{1})|^{2},\text{ if }u=0.
\end{array}
\right.  \label{sym-fft2-abs}%
\end{equation}

It is obvious, that $a(u)$ in (\ref{a333}) can be calculated also as 
\begin{equation}
a(u)=abs(FFT(\textit{fftshift}(J(\tau)))))/N,\text{ } u=1,...,N. \label{a3333}%
\end{equation}

\subsection{Summary}

\begin{enumerate}
\item The all estimates can be calculated as real-valued using cosine
functions as it is in the initial formulations in Proposition 1- Proposition 4 with summation of the observations multiplied by $cos$ with the only difference between these estimates concerning the intervals for
$\tau$, i.e. the intervals of observations.

\item The estimates for $u=0$ are always calculated accordingly to the identical
formulas: see for instance the second lines in (\ref{A_tay_1}),
(\ref{A_tay_12}), (\ref{sym-spectra-est2}).

\item The observation models for $J(\tau)$ are linear equations with respect
to $X(\omega)$ and $|A(\omega)|^{2}$ with a number of equations doubled with
respect to the number of unknown. Thus, the linear algebra methods can be used
in the straightforward manner.

\item The FFT as the fast conventional instrument can be considered as the
preferable algorithm. \ For estimation of the real valued $X(\omega)$ which
can take negative values we should be careful as say in (\ref{a1}) FFT
requires some phase correction and for the symmetric sampling \textit{fftshift} should
be used (\ref{a22}).

\item Applications for the intensity reconstruction $X(\omega)=|A(\omega
)|^{2}$ are simplified as no phase corrections are required. Finally, for the all
three considered cases the reconstructions are identical:%
\begin{align}
a(u)  &  =abs(IFFT(J(\tau))\text{ or }a(u)=abs(FFT(J(\tau)))/N,\label{Final1}%
\\
u  &  =1,...,N,\nonumber
\end{align}
and
\begin{equation}
|A(u)|^{2}=\left\{
\begin{array}
[c]{c}%
a(u+1)\text{, }u=1,2,...,N/2-1\text{,}\\
\\
\frac{1}{4N}\sum_{\tau=-N/2}^{N/2-1}J(\tau)-\frac{1}{2}\sum_{u_{1}=1}%
^{N/2-1}|A(u_{1})|^{2},\text{ if }u=0.
\end{array}
\right.  \label{Final2}%
\end{equation}

The \textit{fftshift} in (\ref{a333}) also can be omitted as it is
equivalent to circular shift which is killed by $abs$. 

As it is shown in
(\ref{a12-final11}), (\ref{a3333}), (\ref{Final1}), IFFT can be replaced for FFT because $J(\tau)$ is even function
of $\tau$.
\end{enumerate}

\section{Holography: setup with reference beam}

\subsection{Non-symmetric sampling}

The discretization of the model (\ref{Integral-Interferogram}) produced as
above according to requirements of FFT, $t=\Delta_{t}\tau$ and $\omega
=\Delta_{\omega}u$, where $\tau$ and $u$ are integers, results in the discrete
observation model:%
\begin{align}
&  J(\tau)=\sum_{u=0}^{N/2-1}|A(u)+R\times\exp(j\frac{2\pi}{N}\tau
u)|^{2},\text{ }\nonumber\\
&  \sum_{u=0}^{N/2-1}(|A(u)|^{2}+|R|^{2}+RA(u)\exp(-j\frac{2\pi}{N}\tau
u)+RA^{\ast}(u)\exp(j\frac{2\pi}{N}\tau u))=\label{Kalenkon_Discrete}\\
&  \sum_{u=0}^{N/2-1}(|A(u)|^{2}+|R|^{2}+2R|A(u)|\cos(\varphi_{A}%
(u)-\frac{2\pi}{N}\tau u)),\nonumber\\
&  \tau=1,..N,\nonumber
\end{align}

where the integer discrete frequency cover the interval $u\in
\{0,1,...,N/2-1\}$, as we assume that the signal $A(u)$ is low band such that $A(u)=0$ for
$u\geq N/2.$

\begin{proposition}
Let the observations $J(\tau)$, $\tau=1,..N$ be defined as in
Eq.(\ref{Kalenkon_Discrete}), then $A(u)$ can be calculated as follows:%
\begin{align}
A(u)  &  =\frac{1}{NR}\sum_{\tau=1}^{N}J(\tau)\exp(j\frac{2\pi}{N}\tau
u)\text{, }0<u\leq N/2-1,\label{Non_SYM_Kalenkov}\\
|A(0)+R|^{2}  &  =\sum_{\tau=1}^{N}J(\tau)-(\sum_{u_{1}=1}^{N/2-1}%
|A(u_{1})|^{2}+\sum_{u_{1}=1}^{N/2-1}R^{2})\text{, }u=0.\nonumber
\end{align}

The calculations with FFT similar to (\ref{A_tay_1_IFFT}) are as follows:
\begin{equation}
a(u)=IFFT(J(\tau))\times\exp(j\frac{2\pi}{N}(u-1)),\text{ }u=1,...,N,
\label{a_kalenkov_nonsym}%
\end{equation}
and%
\begin{align}
A(u)  &  =\frac{1}{R}a(u+1)\text{, if }u=1,..,N/2-1,\label{Kalenkov_IFFT}\\
|A(0)+R|^{2}  &  =\frac{1}{N}\sum_{\tau=1}^{N}J(\tau)-(\sum_{u_{1}=1}%
^{N/2-1}|A(u_{1})|^{2}+\sum_{u_{1}=1}^{N/2-1}R^{2})\text{, if
}u=0\text{.}\nonumber
\end{align}

\end{proposition}
The phase correction in (\ref{a_kalenkov_nonsym}) is as in (\ref{A_tay_1_IFFT}%
) but the operation $real$ is dropped.
The proof of these formulas can be found in Appendix.

If the reference beam $R$ is known, $A(u)$ is defined in the first row of
(\ref{Kalenkov_IFFT}) for $u>0$ and $A(0)$ from the second row. There is no a
unique solution for the complex-valued $A(0)$ from the second equation
(\ref{Kalenkov_IFFT}). However, if we assume that $A(0)$ is real non-negative,
what is natural for the zero frequency of the object, then
the unique solution is of the form
\begin{equation}
A(0)=\sqrt{\frac{1}{N}\sum_{\tau=1}^{N}J(\tau)-(\sum_{u_{1}=1}^{N/2-1}%
|A(u_{1})|^{2}+\sum_{u_{1}=1}^{N/2-1}R^{2})}-R. \label{A_0}%
\end{equation}

Note, that the expression under the square root is always positive.

In what follows, we show simulation tests illustrating the derived estimates
for the complex valued $A(u)$. It is assumed that $A(0)$ is real-valued non-negative.

\begin{figure}[ht]%
\centering
\includegraphics[
height=2.5763in,
width=4.5071in
]%
{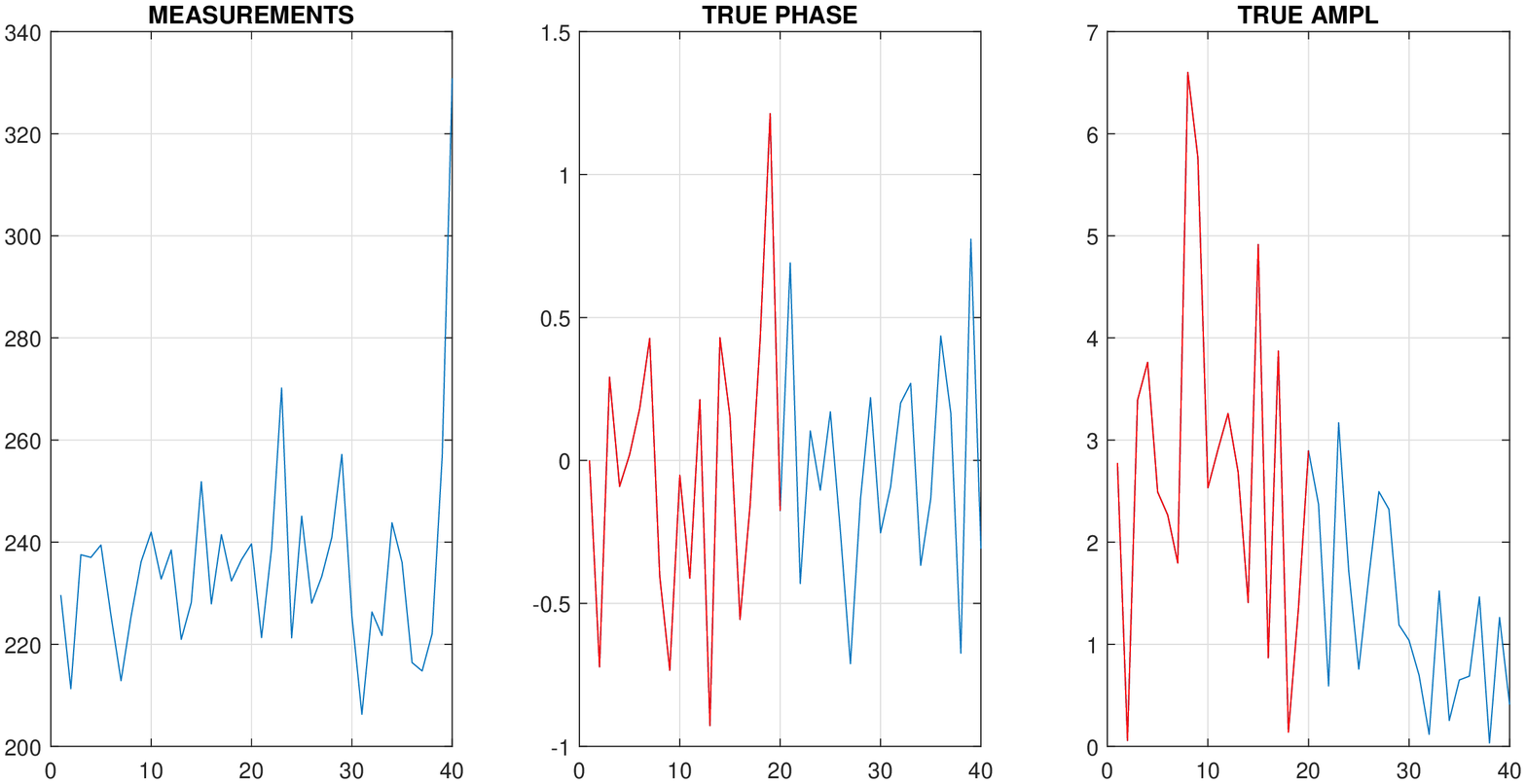}%
\caption{Measurements $J(\tau)$ and true phase and amplitude for the
non-symmetric sampling.}%
\label{spectr_illustr_kalenkov}%
\end{figure}

In Fig.\ref{spectr_illustr_kalenkov} we show the measurements $J(t)$, the
true amplitudes of $A$ (non-negative random values with uniform distribution)
and the true phases of $A$ (random Gaussian with $\sigma=.5).$ The number of
observations $N=40$. According to the above results the perfect reconstruction
is guaranteed provided that $A(u)=0$ for $u\geq N/2$ and for $0\leq u\leq
N/2-1$. The observations are obtained from $A$ are shown in Fig.\ref{spectr_illustr_kalenkov} by zeroing values for $u\geq N/2$. This true
values of the amplitude and phase are marked by red for $0\leq u\leq
N/2-1$.%

\begin{figure}[ht]%
\centering
\includegraphics[
height=2.5763in,
width=4.5071in
]%
{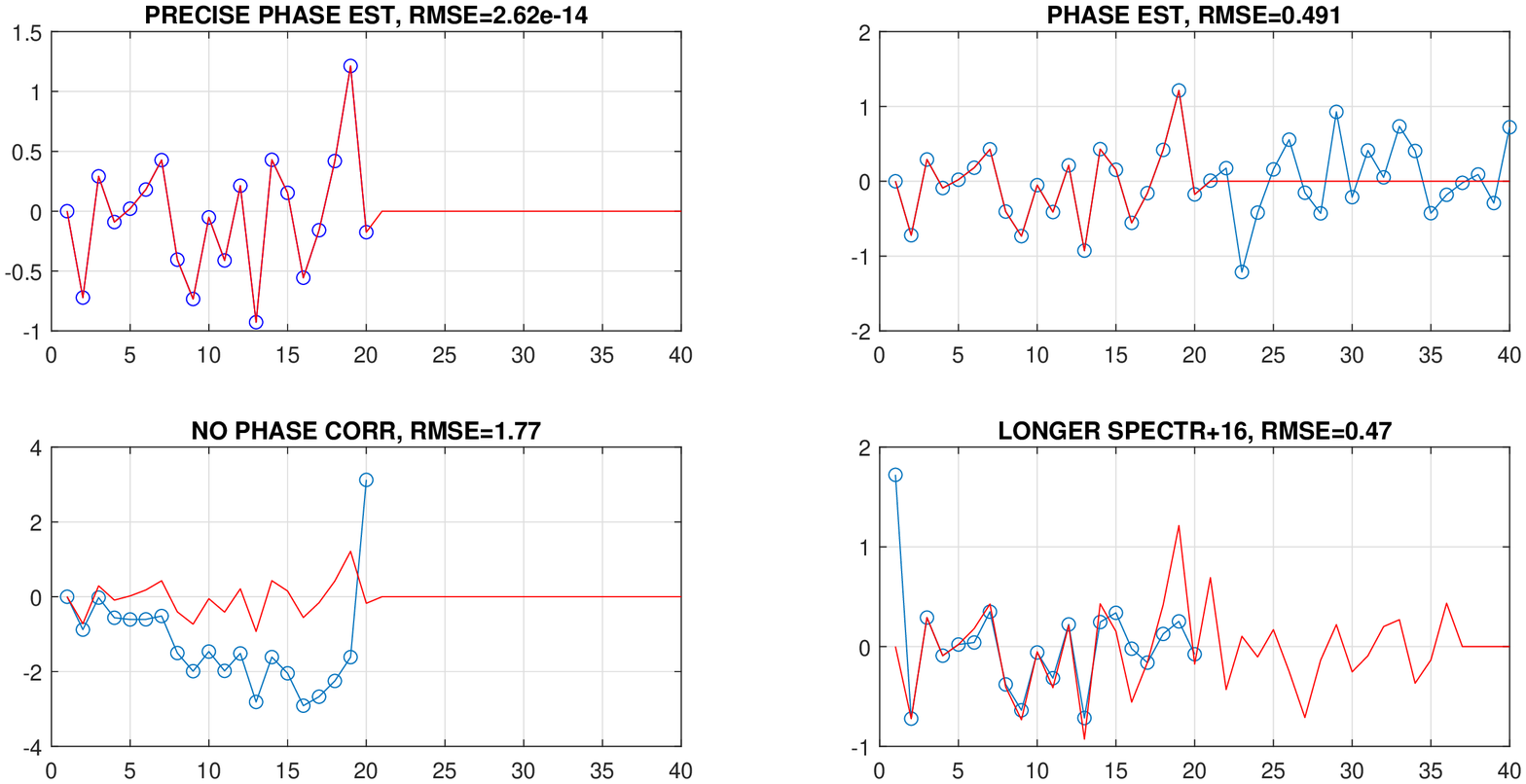}%
\caption{Phase reconstructions: non-symmetric sampling.}%
\label{estimate_illustr_a_phi}%
\end{figure}

\begin{figure}[ht]%
\centering
\includegraphics[
height=2.5763in,
width=4.5071in
]%
{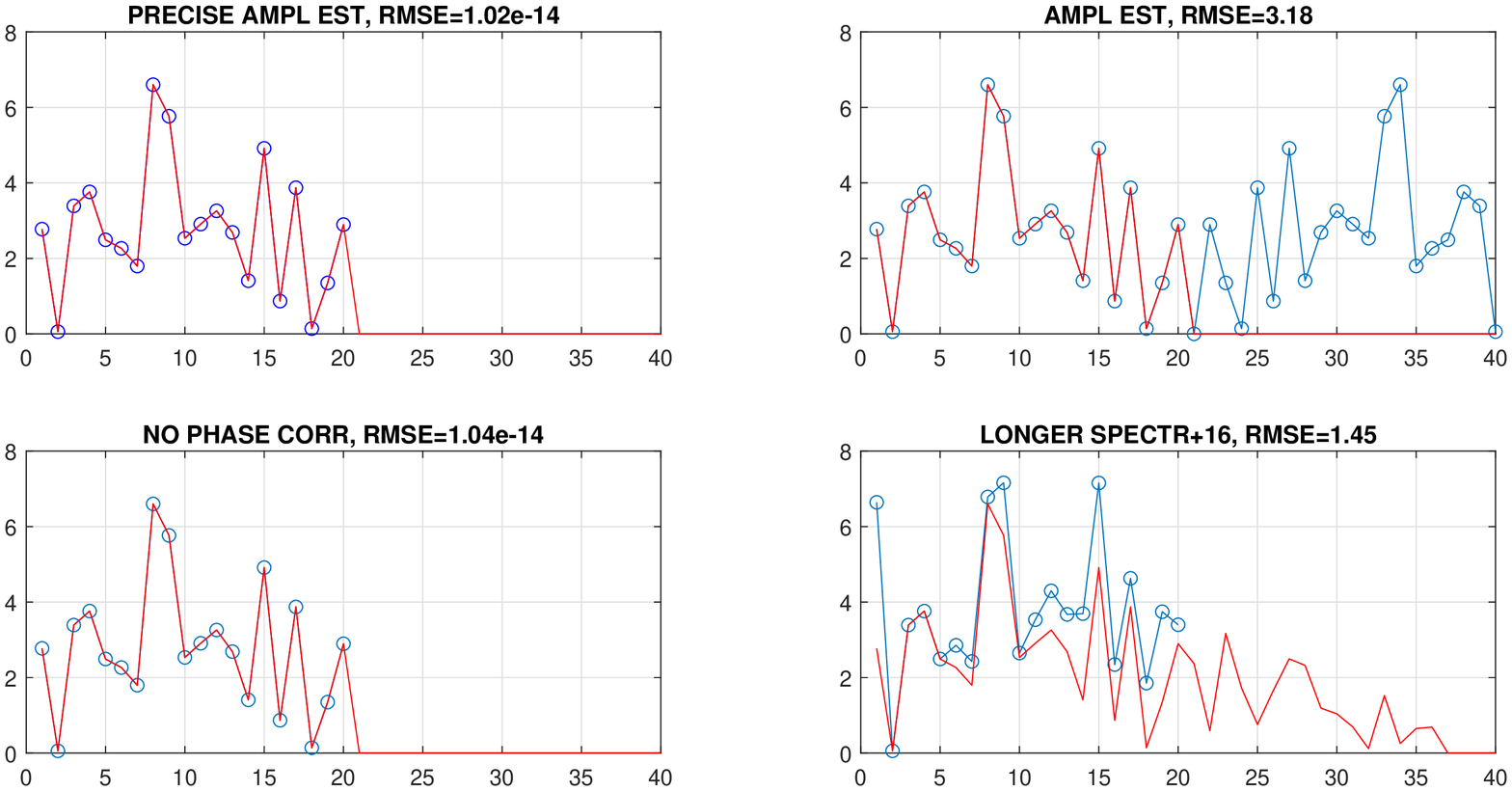}%
\caption{Amplitude reconstructions: non-symmetric sampling.}%
\label{estimate_illustr_a_abs}%
\end{figure}

The precise estimates shown in Figs. \ref{estimate_illustr_a_phi} and
\ref{estimate_illustr_a_abs} are obtained for the truncated spectrum of the
signal which takes zero values for $u\geq N/2$. 

The estimates for the phase
are shown in the upper row of the images in Fig. \ref{estimate_illustr_a_phi}.
In the left image of this row we show the estimates calculated for the basic
interval $0\leq u\leq N/2-1.$ The perfect reconstruction is demonstrated
in this case with the very small value of RMSE. In the right image, we show
the estimates calculated for all $0\leq u\leq N-1$. The red color shows the
true values of the signal. We can see that the estimates calculated for $u\geq
N/2$ are completely erroneous very different from the true values equal to
zero. The RMSE is calculated for the second part of the $\ u$ interval and it
is very large.

There is an obvious symmetry of the estimates for $u\geq N/2$ with respect to
the estimates for $u<N/2$. The right part of the estimate can be interpreted as the negative mirror reflection with respect to the point $N/2+1$ with the phase values
taken with the sign minus.

The second row in Fig. \ref{estimate_illustr_a_phi} illustrates the effects of
two different factors. The left estimate is calculated as it in the first row
but without the phase correction \ for $IFFT$ in Eq.(\ref{a_kalenkov_nonsym}). One
may note that this phase correction is of importance as the estimate is
completely destroyed, what can be noticed from the large values of RMSE as
well as visually comparing the estimate and the true values.%

The last image in Fig. \ref{estimate_illustr_a_phi} shows the influence of the
condition $A(u)=0$ for $u\geq N/2$ for the phase estimation. Here we assume
that $A(u)=0$ for $u\geq N/2+16$. Thus, $16$ pixels of this enlarged $A(u)$
influence the observations $J$, which are naturally different from those in
Fig.\ref{spectr_illustr_kalenkov}. We may note that the phase estimates are
quite sensitive to this larger length of the object spectrum\ $A(u)$ over the
basic interval $0\leq u\leq N/2-1$, at least the corresponding RMSE takes
quite high value.

The estimates for the amplitude shown in Fig. \ref{estimate_illustr_a_abs} are
obtained for the truncated spectrum of the signal which takes zero values for
$u\geq N/2$. In the left image of the upper row we show the estimates
calculated for the basic interval $0\leq u\leq N/2-1.$ The perfect
reconstruction is demonstrated in this case with the very small value of RMSE.
In the right image, we show the estimates calculated for all $0\leq u\leq
N-1$. The red color shows the true values of the signal. We can see that the
estimates calculated for $u\geq N/2$ are completely erroneous and very
different from the true values which are zero for this interval. The RMSE is
calculated for the second part of the $u$ interval and it is very large.

There is an obvious symmetry of the estimates for $u\geq N/2$ with respect to
the estimates for $u<N/2$. The right part of the estimate can be interpreted as the mirror
reflection with respect to the point $N/2+1$ with values taken with the positive
sign what is different with respect to the symmetry of the estimates for
the phase where the negative sign is appeared in reflection. %

The symmetry properties of the estimates discussed here and above follow from
the properties of FFT and IFFT.

The second row in Fig. \ref{estimate_illustr_a_abs} illustrates the effects of
two different factors. The left estimate is calculated as it in the first row
but without the phase correction \ for $IFFT$ in Eq.(\ref{a_kalenkov_nonsym}). One
may note that this phase correction naturally has no influence for the
amplitude estimate.

The last image shows the influence of the non-zero values of $A(u)=0$ for
$u\geq N/2.$ Here we assume that $A(u)=0$ for $u\geq N/2+16$. Thus, $16$
pixels of this enlarged $A(u)$ influence the observations $J$, which are
naturally different from those in Fig.\ref{spectr_illustr_kalenkov}.  We may note that the amplitude estimate is quite
sensitive to this extra data with RMSE taking a high value.

Similar to the case shown in Eq.(\ref{SpectrallyResolvedObs2}) consider a
shifted set of the observations $J(\tau)$ in the form%
\begin{align}
J(\tau)  &  =\sum_{u=0}^{N/2-1}(|A(u)|^{2}+|R|^{2}+2R|A(u)|\cos(\varphi
_{A}(u)-\frac{2\pi}{N}\tau u),\text{ }\label{SpectrallyResolvedObs22}\\
\tau &  =0,1,..N-1,\nonumber
\end{align}
which is different\ from (\ref{Kalenkon_Discrete}) by the set of shifts $\tau$
including now $\tau=0$.

\begin{proposition}
Let the observations $J(\tau)$, $\tau=0,..N-1,$ be defined as in
Eq.(\ref{SpectrallyResolvedObs22}), then $A(u)$ can be calculated as follows:%
\begin{align}
A(u)  &  =\frac{1}{NR}\sum_{\tau=0}^{N-1}J(\tau)\exp(j\frac{2\pi}{N}\tau
u)\text{, }0<u\leq N/2-1,\label{Non_Sym_Kalenkov2}\\
|A(0)+R|^{2}  &  =\sum_{\tau=0}^{N-1}J(\tau)-(\sum_{u_{1}=1}^{N/2-1}%
|A(u_{1})|^{2}+\sum_{u_{1}=1}^{N/2-1}R^{2})\text{, }u=0.\nonumber
\end{align}

\end{proposition}

The calculations with FFT are similar but different from  shown in (\ref{Kalenkov_IFFT}) with
\begin{equation}
a(u)=IFFT(J(\tau)),\text{ }u=1,...,N, \label{a_KALENKOV_2}%
\end{equation}
then%
\begin{align}
A(u)  &  =\frac{1}{R}a(u+1)\text{, if }%
u=1,..,N/2-1,\label{NonSYM-KALENKOV-NO-SHIFT}\\
|A(0)+R|^{2}  &  =\frac{1}{N}\sum_{\tau=1}^{N}J(\tau)-(\sum_{u_{1}=1}%
^{N/2-1}|A(u_{1})|^{2}+\sum_{u_{1}=1}^{N/2-1}R^{2})\text{,
}u=0\text{.}\nonumber
\end{align}

Comparing (\ref{a_kalenkov_nonsym}) with (\ref{a_KALENKOV_2}) we can recognize
the main difference between these two estimates as the calculation of $a(u)$
for (\ref{NonSYM-KALENKOV-NO-SHIFT}) does not requite the phase correction
essential in (\ref{a_kalenkov_nonsym}).

We do not provide the proof of this proposition as it is similar to shown in
Appendix for Proposition 5.

The numerical illustrations for this algorithm are also similar to shown in
Figs. \ref{spectr_illustr_kalenkov} - \ref{estimate_illustr_a_abs}, naturally,
without the phase correction analysis which is not required in
(\ref{a_KALENKOV_2}).

\subsection{\textbf{Symmetric sampling }}

The measurement equations are
\begin{align}
&  J(\tau)=\sum_{u=0}^{N/2-1}|A(u)+R\exp(j\frac{2\pi}{N}\tau u)|^{2},\text{
}\label{sym-sampling-reference}\\
&  \sum_{u=0}^{N/2-1}(|A(u)|^{2}+|R|^{2}+RA(u)\exp(-j\frac{2\pi}{N}\tau
u)+RA^{\ast}(u)\exp(j\frac{2\pi}{N}\tau u)),\nonumber\\
&  \tau=-N/2,..N/2-1.\nonumber
\end{align}

\begin{proposition}
Let the observations $J(\tau)$ be defined as in (\ref{sym-sampling-reference}%
), then $A(u)$ can be calculated as follows:
\begin{equation}%
\begin{array}
[c]{c}%
A(u)=\frac{1}{RN}\sum_{\tau=-N/2}^{N/2-1}J(\tau)\exp(j\frac{2\pi}{N}\tau
u)\text{, if }u=1,...,N/2-1,\\
\\
|A(0)+R|^{2}=\frac{1}{N}\sum_{\tau=-N/2}^{N/2-1}J(\tau)-(\sum_{u_{1}%
=1}^{N/2-1}|A(u_{1})|^{2}+\sum_{u_{1}=1}^{N/2-1}R^{2})\text{, }u=0.
\end{array}
\label{REC-SYMM-REFERENCE}%
\end{equation}
If FFT is used for calculations these estimates can be presented in the following two forms:

\begin{equation}
a(u)=IFFT(J(\tau))\times\exp(-j\pi(u-1)), u=1,...,N,
\label{x_u_old_old}%
\end{equation}

or
\begin{equation}
a(u)=IFFT(\textit{fftshift}(J(\tau))),\text{ }u=1,...,N, \label{est_FFTSHIFT_with_shift}\\%
\end{equation}.

then, 
\begin{align}
A(u)  &  =\frac{1}{R}a(u+1)\text{, }u=1,2,...,N/2-1\text{,}\label{est_FFTSHIFT}\\
& \nonumber\\
|A(0)+R|^{2}  &  =\frac{1}{N}\sum_{\tau=-N/2}^{N/2-1}J(\tau)-(\sum_{u_{1}%
=1}^{N/2-1}|A(u_{1})|^{2}+\sum_{u_{1}=1}^{N/2-1}R^{2})\text{,
}u=0.\nonumber
\end{align}

\end{proposition}
 The proof is given in Appendix.

In our simulation experiments with this latter algorithm  (\ref{est_FFTSHIFT_with_shift})-(\ref{est_FFTSHIFT})
using \textit{fftshift,} we assume that $A(0)
$ is real positive, what allows to estimate $A(0)$ as in (\ref{A_0}).

\begin{figure}[ht]%
\centering
\includegraphics[
height=3.0865in,
width=4.6037in
]%
{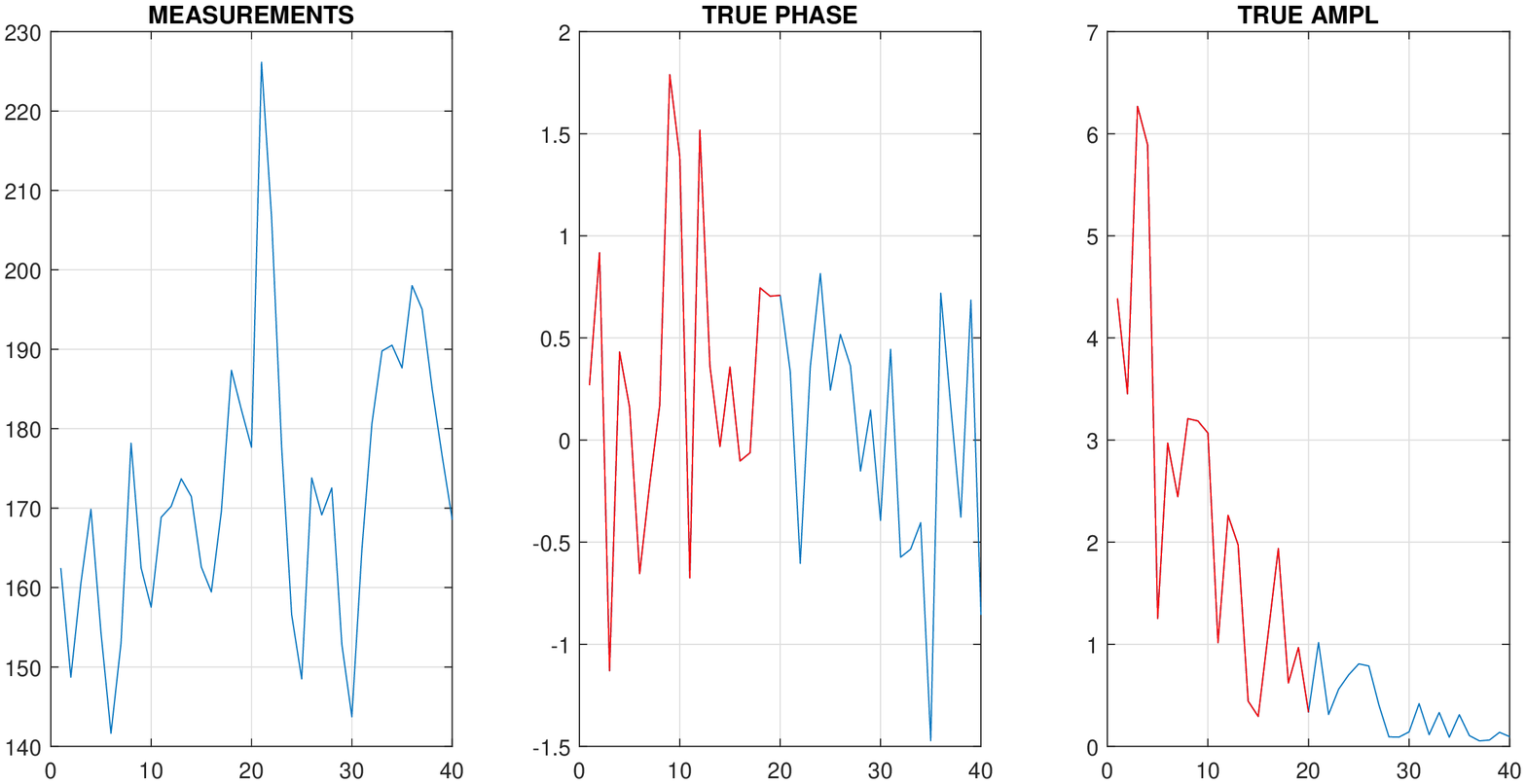}%
\caption{Observations and true \ phase/amplitude for symmetric sampling.}%
\label{spectr_symm_kalenkov}%
\end{figure}

The observations and the true amplitude/phase images are shown in
Fig.\ref{spectr_symm_kalenkov}.

\begin{figure}[ht]%
\centering
\includegraphics[
height=3.0865in,
width=4.6037in
]%
{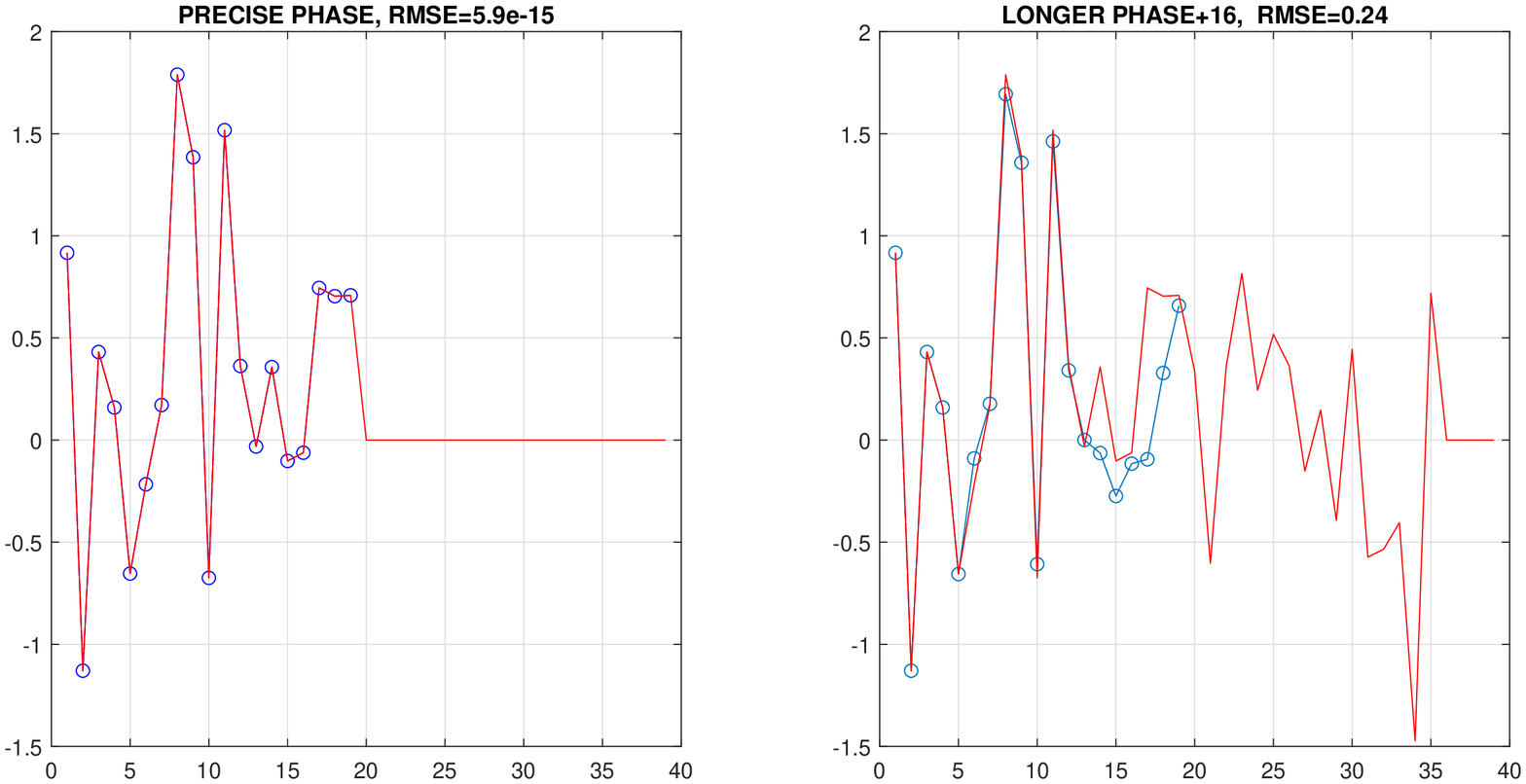}%
\caption{Phase reconstructions, effect of the longer nonzero object spectrum.}%
\label{estimate_symm_kalenkov_phase}%
\end{figure}

\begin{figure}[ht]%
\centering
\includegraphics[
height=3.0865in,
width=4.6037in
]%
{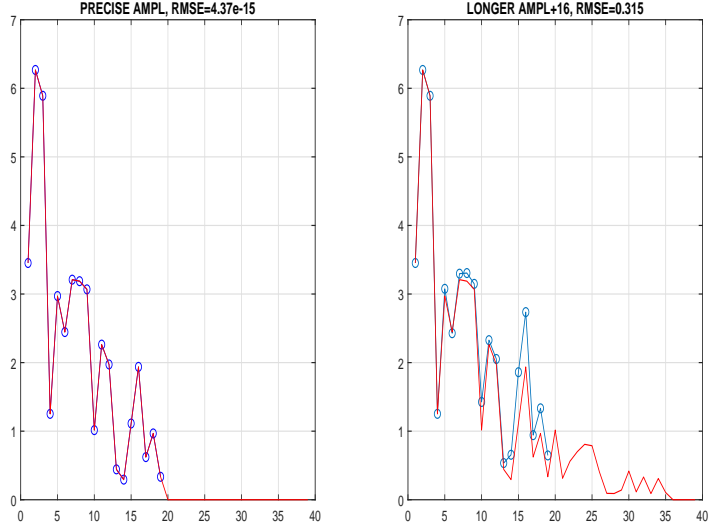}%
\caption{Amplitude reconstructions, effect of the longer nonzero object
spectrum.}%
\label{estimate_symm_kalenkov_amplit}%
\end{figure}

Illustrations of the sensitivity of the estimates to $|A(u_{1})|^{2}>0$ for
$u>N/2$ are given in Figs.\ref{estimate_symm_kalenkov_phase} and
Figs.\ref{estimate_symm_kalenkov_amplit} for phase and amplitude
reconstructions, respectively. In the left images, the results demonstrating the
perfect reconstruction are shown for the lower band spectrum and for the
degradation of the reconstruction are in the right images where the spectrum
is wider than $N/2$ by $16$ pixels. In the latter case visually and
numerically the degradation of the estimates is essential.%

\begin{figure}[ht]%
\centering
\includegraphics[
height=3.0865in,
width=4.6037in
]%
{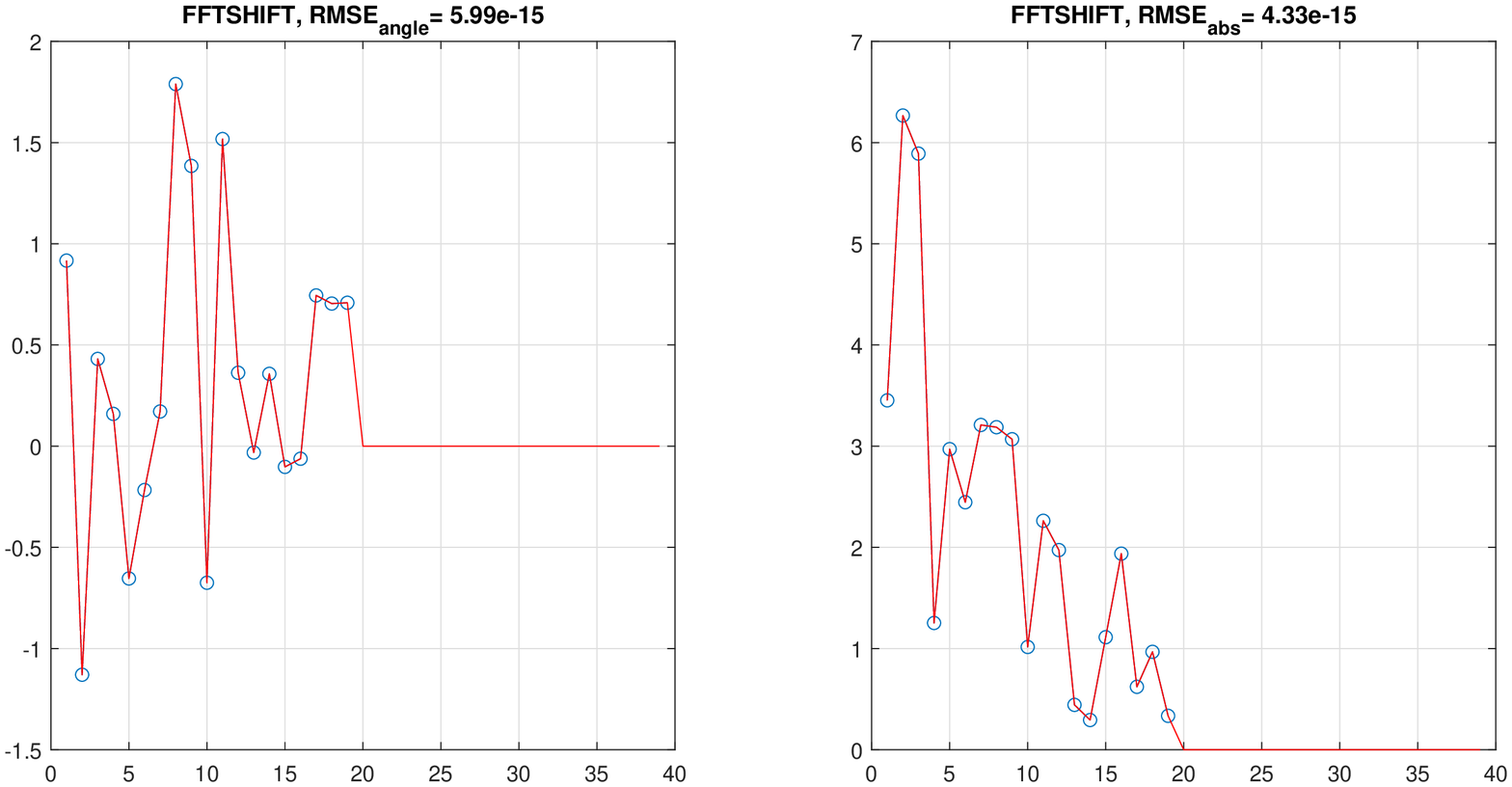}%
\caption{Perfect phase and amplitude reconstructions by IFFT with the \textit{fftshift}
operation.}%
\label{estimate_symm_kalenkov_fftshif}%
\end{figure}

The estimates in Fig.\ref{estimate_symm_kalenkov_phase} and
Fig.\ref{estimate_symm_kalenkov_amplit} are produced using the algorithm (\ref{x_u_old_old}) with the phase correction for IFFT. In Fig.\ref{estimate_symm_kalenkov_fftshif} we show the results obtained by the algorithm (\ref{est_FFTSHIFT_with_shift}) with the \textit{fftshift} operation. The comparison with Fig.\ref{estimate_symm_kalenkov_phase} and
Fig.\ref{estimate_symm_kalenkov_amplit} confirms the identity of the results delivered by these algorithms.


\subsection{Summary}

\begin{enumerate}
\item The all estimates can be calculated as FFT or IFFT of the observations, see
(\ref{a_kalenkov_nonsym}), (\ref{a_KALENKOV_2}), (\ref{est_FFTSHIFT_with_shift}%
) with the summation interval defined by the observation shifts.

\item The estimates for $u=0$ are always calculated accordingly to the similar
formulas: see for instance the second lines in (\ref{Non_SYM_Kalenkov}),
(\ref{Non_Sym_Kalenkov2}), (\ref{REC-SYMM-REFERENCE}).

\item The FFT as the fast conventional instrument can be considered as a
preferable algorithm. The most compact form of estimation is achieved using
\textit{fftshift} together with $IFFT$ in (\ref{est_FFTSHIFT_with_shift}).
\item  The MATLAB demo-codes reproducing the presented numerical tests are publicly available: http://www.cs.tut.fi/sgn/imaging/sparse/.
\end{enumerate}

\section{\textbf{Appendix}}

\subsection{\textit{The proof of Proposition 1}}

Remind, that the Kronecker delta-function defined for integer $n,m$ as
\begin{equation}
\frac{1}{N}\sum_{k=1}^{N}\exp(j2\pi(n-m)k/N)=\delta_{n,m} \label{Kronecker1}%
\end{equation}
is periodic on $(n-m)$. We apply it only for the\ single period
\begin{equation}
|n-m|<N. \label{Kronecker_Cond}%
\end{equation}
Using this definition of the Kronecker function and the expression%
\begin{align}
\cos(\frac{2\pi}{N}\tau u)\cos(\frac{2\pi}{N}\tau v)  &  =\frac{1}{4}%
[(\exp(j\frac{2\pi}{N}\tau u)+\exp(-j\frac{2\pi}{N}\tau u))\label{cosine}\\
&  +(\exp(j\frac{2\pi}{N}\tau v)+\exp(-j\frac{2\pi}{N}\tau v))],\nonumber
\end{align}
simple calculations show that for integer $u$ and $v$%
\begin{align}
&  \frac{2}{N}\sum_{\tau=1}^{N}\cos(\frac{2\pi}{N}\tau u)\cos(\frac{2\pi}%
{N}\tau v)=\delta_{u,v}+\delta_{u,-v},\label{Kronecker}\\
&  u,v\in\lbrack0,..,N/2-1].\nonumber
\end{align}
The later restrictions on $u$ and $v$ are essential as they guarantee that
$|u-v|<N$ and $|u+v|<N$ and only the single period of the Kronecker function
is in use for (\ref{cosine}).

Let us calculate $\dfrac{1}{N}\sum_{\tau=1}^{N}J(\tau)\cos(\frac{2\pi}{N}\tau
u)$ inserting the expression for $J(\tau)$, Eq.(\ref{SpectrallyResolvedObs}), and
using Eq.(\ref{Kronecker}).

Then,%
\begin{align}
&  \frac{2}{N}\sum_{\tau=1}^{N}\sum_{u_{1}=0}^{N/2-1}X(u_{1})(1+\cos
(\frac{2\pi}{N}\tau u_{1}))\cos(\frac{2\pi}{N}\tau u)=\nonumber\\
&  \sum_{u_{1}=0}^{N/2-1}X(u_{1})[\frac{2}{N}\sum_{\tau=1}^{N}\cos(\frac{2\pi
}{N}\tau u)+\frac{2}{N}\sum_{\tau=1}^{N}\cos(\frac{2\pi}{N}\tau u_{1}%
)\cos(\frac{2\pi}{N}\tau u)]=\nonumber\\
&  \sum_{u_{1}=0}^{N/2-1}X(u_{1})[2\delta_{u,0}+\delta_{u,u_{1}}%
+\delta_{u,-u_{1}}]=2\sum_{u_{1}=0}^{N/2-1}X(u_{1})\delta_{u,0}+|A(u)|^{2}%
+|A(0)|^{2}\delta_{u,0}. \label{proof0}%
\end{align}
Thus, we obtain that
\begin{equation}
\frac{1}{N}\sum_{\tau=1}^{N}J(\tau)\cos(\frac{2\pi}{N}\tau u)=2\sum_{u_{1}%
=0}^{N/2-1}X(u_{1})\delta_{u,0}+X(u)+X(0)\delta_{u,0}. \label{proof}%
\end{equation}
The first line of (\ref{A_tay_1}) follows from (\ref{proof}) if $u=1,..,N/2-1
$. For $u=0,$ we obtain from (\ref{proof})%
\begin{align}
&  \frac{1}{N}\sum_{\tau=1}^{N}J(\tau)=2\sum_{u_{1}=0}^{N/2-1}X(u_{1}%
)+X(0)+X(0)=\label{inter}\\
&  4X(0)+2\sum_{u_{1}=1}^{N/2-1}X(u_{1}).\nonumber
\end{align}
and the second line of (\ref{A_tay_1}) follows from the last formulas.

We can use for (\ref{A_tay_1}) the complex exponent, then%
\begin{equation}
X(u)=\left\{
\begin{array}
[c]{c}%
\frac{1}{N}\operatorname{real}(\sum_{\tau=1}^{N}J(\tau)\exp(+j\frac{2\pi}%
{N}\tau u))\text{, }1\leq u\leq N/2-1,\\
\\
\frac{1}{4N}\sum_{\tau=1}^{N}J(\tau)-\frac{1}{2}\sum_{u_{1}=1}^{N/2-1}%
X(u_{1})\text{, }u=0,
\end{array}
\right.  \label{Complex}%
\end{equation}
then the calculations can be produced using FFT as follows%
\begin{align*}
x(u)  &  =\operatorname{real}(IFFT(J(\tau))\times\exp(j\frac{2\pi}%
{N}(u-1))),\\
u  &  =1,..,N.
\end{align*}
Then, the formula (\ref{Complex}) can
be rewritten as
\[
X(u)=\left\{
\begin{array}
[c]{c}%
x(u+1)\text{, if }u=1,..,N/2-1,\\
\\
\dfrac{1}{4N}\sum_{\tau=1}^{N}J(\tau)-\frac{1}{2}\sum_{u_{1}=1}^{N/2-1}%
X(u_{1}),\text{ if }u=0.
\end{array}
\right.
\]
what confirms (\ref{A_tay_1_IFFT}).

It completes the proof of Proposition 1.

\subsection{\textit{The proof of Proposition 4.}}

For the symmetric sampling observation model is of the form:
\begin{align}
&  J(\tau)=\sum_{u=0}^{N/2-1}X(u)(2+\exp(-j\frac{2\pi}{N}\tau u)+\exp
(j\frac{2\pi}{N}\tau u)),\label{Prop3_obs}\\
&  \tau=-N/2,...,N/2-1.\nonumber
\end{align}
We calculate DFT of $J(\tau)$ as follows
\begin{align}
&  \frac{1}{N}\sum_{\tau=-N/2}^{N/2-1}J(\tau)\exp(j\frac{2\pi}{N}\tau
u_{1})=\label{Prop3_Proof}\\
&  \sum_{u=0}^{N/2-1}\frac{1}{N}\sum_{\tau=-N/2}^{N/2-1}X(u)(2+\exp
(-j\frac{2\pi}{N}\tau u)+\exp(j\frac{2\pi}{N}\tau u))\exp(j\frac{2\pi}{N}\tau
u_{1})=\nonumber\\
&  \sum_{u=0}^{N/2-1}[X(u)(2\delta_{u_{1},0}+\delta_{u,u_{1}}+\delta
_{u,-u_{1}}]=\nonumber\\
&  \sum_{u=0}^{N/2-1}X(u)2\delta_{u_{1},0}+X(u_{1})^{2}+X(0)\delta_{u_{1}%
,0}.\nonumber
\end{align}
It follows that%
\begin{align}
X(u_{1})  &  =\frac{1}{N}\operatorname{real}(\sum_{\tau=-N/2}^{N/2-1}%
J(\tau)\exp(j\frac{2\pi}{N}\tau u_{1}))\text{, }1\leq u_{1}\leq
N/2-1,\label{Prop3_results}\\
& \nonumber\\
X(0)  &  =(\frac{1}{N}\sum_{\tau=-N/2}^{N/2-1}J(\tau)-2\sum_{u=1}%
^{N/2-1}|A(u)|^{2})/4,\text{ }u_{1}=0.\nonumber
\end{align}
The last equation follows from Eq.(\ref{Prop3_Proof}) for $u_{1}=0$. It proves
Eqs.(\ref{sym-spectra-est2}).

In order to derive the formulas with FFT, we produce a change of the summation
variable in the first line of (\ref{Prop3_results}):%
\begin{align}
&  \frac{1}{N}\sum_{\tau=-N/2}^{N/2-1}J(\tau)\exp(j\frac{2\pi}{N}\tau
u_{1})=\frac{1}{N}\sum_{\tau_{1}=0}^{N-1}J(\tau_{1}-N/2)\exp(j\frac{2\pi}%
{N}(\tau_{1}-N/2)u_{1})=\label{replacement}\\
&  \frac{1}{N}\sum_{\tau_{1}=0}^{N-1}J(\tau_{1}-N/2)\exp(j\frac{2\pi}{N}%
(\tau_{1}u_{1}))\exp(-j\pi u_{1})=\nonumber\\
&  IFFT(J(\tau))\times\exp(-j\pi(u_{1}-1)),\text{ }1\leq u_{1}\leq
N/2-1.\nonumber
\end{align}
where%
\[
IFFT(J(\tau))=\frac{1}{N}\sum_{\tau=1}^{N}J(\tau)\exp(j\frac{2\pi}{N}%
(\tau-1)(u-1)).
\]
The last expression in (\ref{replacement}) proves Eq. (\ref{x_u_old}).

In order to prove (\ref{a22})-(\ref{sym-fft2}) note, that the sequence
$J(\tau)$ in (\ref{sym-spectra-est2}) is ordered according to $\tau\in
\lbrack-N/2,-N/2+1,...,-1,0,1,...,N/2-1]$. The MATLAB \textit{fftshift
}operation in (\ref{a22}) circularly shifts $\tau=0$ to the initial position
and the sequence takes the form $\tau\in\lbrack
0,1,...,N/2-1,-N/2,-N/2+1,...,-1]$. The $J(\tau)$ is periodic with the period
equal to $N$, $J(\tau)=J(\tau+N).$ We add $N$ to the negative items of the
obtained sequence and arrive to another sequence $\tau\in\lbrack0,1,..,N-1]$
with the values of $J(\tau)$ identical to those for the initial sequence
$\tau\in\lbrack-N/2,-N/2+1,...,-1,0,1,...,N/2-1]$.

In order to extract the values $X(u)$ from this sequence we use the results of
Proposition 2 dealing exactly with such sequence of $J(\tau)$, $\tau\in
\lbrack0,1,..,N-1]$. Then, the formulas (\ref{a12})-(\ref{A_tay_1_IFFT2}) can
be used for this reconstruction. 

It proves Eqs. (\ref{a22})-(\ref{sym-fft2}) and
completes the proof of the proposition.

\subsection{\textit{The proof of Proposition 5. }}

Observations are of the form%
\begin{align}
&  J(\tau)=\sum_{u=0}^{N/2-1}|A(u)+R\exp(j\frac{2\pi}{N}\tau u)|^{2},\text{
}\tau=1,..N.\label{J_t}\\
&  \sum_{u=0}^{N/2-1}(|A(u)|^{2}+R^{2}+RA(u)%
\operatorname{A}%
\exp(-j\frac{2\pi}{N}\tau u)+RA^{\ast}(u)\exp(j\frac{2\pi}{N}\tau
u)).\nonumber
\end{align}
Let us calculate DFT for $J(\tau):$%
\begin{align*}
&  \frac{1}{N}\sum_{\tau=1}^{N}J(\tau)\exp(j\frac{2\pi}{N}\tau u_{1})=\\
&  \sum_{u=0}^{N/2-1}\frac{1}{N}\sum_{\tau=1}^{N}[(|A(u)|^{2}+R^{2}%
+RA(u)\exp(-j\frac{2\pi}{N}\tau u)+\\
&  RA^{\ast}(u)\exp(j\frac{2\pi}{N}\tau u))]\exp(j\frac{2\pi}{N}\tau
u_{1})=\\
&  \sum_{u=0}^{N/2-1}[(|A(u)|^{2}\delta_{u_{1},0}+R^{2}\delta_{u_{1}%
,0}+RA(u)\delta_{u,u_{1}}+RA^{\ast}(u)\delta_{u,-u_{1}}]=\\
&  \sum_{u=0}^{N/2-1}|A(u)|^{2}\delta_{u_{1},0}+\sum_{u=0}^{N/2-1}%
R^{2}\delta_{u_{1},0}+RA(u_{1})+RA^{\ast}(0)\delta_{u_{1},0}].
\end{align*}
Here we use the definition of the Kronecker $\delta$-function
(\ref{Kronecker1}). It results in
\begin{align}
RA(u_{1})  &  =\frac{1}{N}\sum_{\tau=1}^{N}J(\tau)\exp
(j\frac{2\pi}{N}\tau u_{1})\text{, }1\leq u_{1}\leq
N/2-1,\label{NON_SYM_Lalenkov_Proof}\\
|A(0)+R(0)|^{2}  &  =\frac{1}{N}\sum_{\tau=1}^{N}J(\tau)-(\sum_{u=1}%
^{N/2-1}|A(u)|^{2}+\sum_{u=1}^{N/2-1}R^{2}),\text{ }u_{1}=0.\nonumber
\end{align}
Using $IFFT$ we arrive to the formula
\begin{align*}
A(u)  &  =\frac{1}{R}IFFT(J(\tau))\exp(j\frac{2\pi}{N}\tau (u-1))\text{, }1\leq
u\leq N/2-1,\\
|A(0)+R(0)|^{2}  &  =\frac{1}{N}\sum_{\tau=1}^{N}J(\tau)-(\sum_{u_{1}%
=1}^{N/2-1}|A(u_{1})|^{2}+\sum_{u_{1}=1}^{N/2-1}R^{2})\text{, }u=0,
\end{align*}
what can be written in the form (\ref{Kalenkov_IFFT}).

It completed the proof of Proposition 5.

\subsection{\textit{The proof of Proposition 7}}

Observations are of the form:%
\begin{align*}
&  J(\tau)=\sum_{u=0}^{N/2-1}|A(u)+R\cdot\exp(j\frac{2\pi}{N}\tau
u)|^{2},\text{ }\tau=-N/2,..,N/2-1.\\
&  \sum_{u=0}^{N/2-1}(|A(u)|^{2}+R^{2}+RA(u)\exp(-j\frac{2\pi}{N}\tau
u)+RA^{\ast}(u)\exp(j\frac{2\pi}{N}\tau u)).
\end{align*}
Calculate DFT for these observations%
\begin{align}
&  \frac{1}{N}\sum_{\tau=-N/2}^{N/2-1}J(\tau)\exp(j\frac{2\pi}{N}\tau
u_{1})=\label{Prof4_proof}\\
&  \sum_{u=0}^{N/2-1}\frac{1}{N}\sum_{\tau=-N/2}^{N/2-1}[(|A(u)|^{2}%
+R^{2}+RA(u)\exp(-j\frac{2\pi}{N}\tau u)+\nonumber\\
&  RA^{\ast}(u)\exp(j\frac{2\pi}{N}\tau u))]\exp(j\frac{2\pi}{N}\tau
u_{1})=\nonumber\\
&  \sum_{u=0}^{N/2-1}[(|A(u)|^{2}\delta_{u_{1},0}+R^{2}\delta_{u_{1}%
,0}+RA(u)\delta_{u,u_{1}}+RA^{\ast}(u)\delta_{u,-u_{1}}]=\nonumber\\
&  \sum_{u=0}^{N/2-1}|A(u)|^{2}\delta_{u_{1},0}+\sum_{u=0}^{N/2-1}%
R^{2}\delta_{u_{1},0}+RA(u_{1})+RA^{\ast}(0)\delta_{u_{1}%
,0}].\nonumber
\end{align}
The results are as follows%
\begin{align}
RA(u)  &  =\frac{1}{N}\sum_{\tau=-N/2}^{N/2-1}J(\tau)\exp(j\frac{2\pi}%
{N}\tau u)\text{, }1\leq u\leq N/2-1,\label{Sym_Kalenkov_Proof}\\
|A(0)+R(0)|^{2}  &  =\frac{1}{N}\sum_{\tau=-N/2}^{N/2-1}J(\tau)-(\sum
_{u_{1}=1}^{N/2-1}|A(u_{1})|^{2}+\sum_{u_{1}=1}^{N/2-1}R^{2}),\text{
}u=0.\nonumber
\end{align}
It proves (\ref{REC-SYMM-REFERENCE}).


Let us produce the transformations in the first line equation in Eq.(\ref{Sym_Kalenkov_Proof}) similar to used in the proof of Proposition 4: $\frac{1}{N}\sum_{\tau=-N/2}^{N/2-1}J(\tau)\exp(j\frac{2\pi
}{N}\tau u_{1})=
\frac{1}{N}\sum_{\tau_{1}=0}^{N-1}J(\tau_{1}-N/2)\exp
(j\frac{2\pi}{N}(\tau_{1}u_{1}))\exp(-j\pi u_{1})$ and introduce
\begin{equation}
a(u)=IFFT(J(\tau))\times\exp(-j\pi(u-1)),u\in\lbrack1,...,N], \label{a44}%
\end{equation}
then Eq.(\ref{Sym_Kalenkov_Proof}) takes the form:
\begin{align}
A(u)  &  =\frac{1}{R}a(u+1)\text{, }u=1,2,...,N/2-1\text{,}\label{sym-FFT}\\
& \nonumber\\
|A(0)+R|^{2}  &  =\frac{1}{N}\sum_{\tau=-N/2}^{N/2-1}J(\tau)-(\sum
_{u_{1}=1}^{N/2-1}|A(u_{1})|^{2}+\sum_{u_{1}=1}^{N/2-1}R^{2})\text{,
}u=0.\nonumber
\end{align}
It proves the algorithm with $a(u)$ defined as in (\ref{x_u_old_old}). 

The phase correction by the complex exponent $\exp(-j\pi(u-1))$ is essential
for this form of the solution.

In order to simplify this solution, we involve \textit{fftshift} operation as it is
in the proof of Proposition 4 and arrive to the following compact expression:
\begin{equation}
a(u)=IFFT(\textit{fftshift}(J(\tau))), u\in\lbrack1,...,N], \label{a55}%
\end{equation}
then
\begin{align}
A(u)  &  =\frac{1}{R}a(u+1)\text{, }u=1,2,...,N/2-1\text{,}\label{kalenkov-Proof}\\
& \nonumber\\
|A(0)+R|^{2}  &  =\frac{1}{N}\sum_{\tau=-N/2}^{N/2-1}J(\tau)-(\sum
_{u_{1}=1}^{N/2-1}|A(u_{1})|^{2}+\sum_{u_{1}=1}^{N/2-1}R^{2})\text{,
}u=0.\nonumber
\end{align}
It confirms Eqs.(\ref{est_FFTSHIFT_with_shift})-(\ref{est_FFTSHIFT}) and completes the proof of Proposition 7.

\section{ACKNOWLEDGMENTS}

This work is done as a part of CIWIL project funded by Jane and Aatos Erkko
Foundation and supported by Finnish Academy of Science, Project No. 287150, 2015-2019.

\end{document}